\documentclass{amsart}

\newtheorem{theorem}{Theorem}[section]
\newtheorem{lemma}[theorem]{Lemma}

\newtheorem{corollary}[theorem]{Corollary}
\newtheorem{question}[theorem]{Question}
\newtheorem{problem}[theorem]{Problem}
\newtheorem{open problem}[theorem]{Open Problem}
\newtheorem{conjecture}[theorem]{Conjecture}

\theoremstyle{definition}
\newtheorem{definition}[theorem]{Definition}

\usepackage{amscd,amssymb}
\usepackage{graphicx}

\begin{document}

\title[Finite basis problem]
{Finite basis problem\\ for varieties of algebraic systems}

\author[Vesselin Drensky]
{Vesselin Drensky}
\address{Institute of Mathematics and Informatics,
Bulgarian Academy of Sciences,
1113 Sofia, Bulgaria}
\email{drensky@math.bas.bg}


\maketitle
\begin{abstract}
This is a survey on the finite basis problem for varieties of algebraic systems.
Our exposition is in two directions:

(i) We give numerous examples of varieties which are not finitely based.

(ii) We give examples of important varieties with the property that they and their subvarieties are finitely based.

A special attention is paid on the varieties of semigroups, groups, associative, Lie and other nonassociative algebras.
\end{abstract}

\section{Introduction}

Naively stated, one of the main problems in abstract algebra is:

\begin{problem}\label{naive main problem}
Describe all groups, semigroups, rings and algebras over a field (both associative and nonassociative).
\end{problem}

Clearly, it is hopeless to list all groups, semigroups, rings and algebras and we need another approach.
We shall discuss two possibilities:

$\bullet$ Structure approach.

$\bullet$ Combinatorial approach.

\subsection{Structure theory of algebraic systems}
Figuratively speaking, structure theory is like anatomy.
(The word ``anatomy'' comes from the Greek $\alpha\nu\alpha\tau o\mu\acute{\eta}$ which means ``dissection''.
One studies the structure of the body, its organs and the interactions between the organs.)
A typical example is the theory of finite-dimensional algebras.
If $R$ is an associative algebra over a field $K$ and $\dim(R)<\infty$,
then one defines an ideal $J(R)$ of $R$ called the Jacobson radical such that:

$\bullet$ The factor algebra $R/J(R)$ is semisimple. Every semisimple algebra is a direct sum of simple algebras.
The only finite-dimensional simple algebras are the matrix algebras with entries from a finite-dimensional division ring
(division ring = not necessarily commutative field).

$\bullet$ The radical is nilpotent, i.e. there is a positive integer such that $u_1\cdots u_n=0$ for all $u_1,\ldots,u_n$ in $J$
which can be written as $J^n=0$.

$\bullet$ Over a perfect field $K$ the algebra $R$ is a split extension of a semisimple subalgebra by the radical, i.e.
$R\cong R/J(R)\oplus J(R)$.

\subsection{Combinatorial approach to the theory of algebraic systems}
The combinatorial approach is similar to the Organism Classification of Carl Linnaeus.
The living creatures are divided in classes (e.g. reptilia, birds, mammals) and each class has subclasses.
In the same way the algebraic objects are collected in classes called varieties depending on the identical relations which they satisfy.
The advantage of this approach is that we can treat in the same way groups, semigroups, rings and algebras (commutative, associative, Lie, Jordan, etc.).
This idea comes from the works by Garrett Birkhoff on universal algebras \cite{Birk1} and Bernhard Neumann on groups \cite{Neu}.

\subsection{Varieties of algebraic systems}
Recall the definitions of an $\Omega$-algebra and the main objects in our considerations.

\begin{definition}\label{definition of Omega-algebras}
{\rm (i)} Let $\Omega$ be a set of operations
which is a union of subsets of $n$-ary operations $\Omega_n$, $n=0,1,2,\ldots$, subject to a number of identities.
A set $A$ is an $\Omega$-algebra (or an algebraic system of signature $\Omega$) if for every $\omega\in\Omega_n$, $n=0,1,2,\ldots$, we have a mapping
\[
\omega:\underbrace{A\times\cdots\times A}_{n\;\rm{times}}\to A
\]
satisfying the identities in the definition of $\Omega$.

{\rm (ii)} Let $\mathfrak V$ be the class of all $\Omega$-algebras.
For any set $X$ there is an universal $\Omega$-algebra $F(X)$ (the free $\Omega$-algebra generated by $X$) such that for any $A\in\mathfrak V$ any mapping $X\to A$
can be extended in a unique way to a homomorphism of $\Omega$-algebras $F(X)\to A$.
\end{definition}

For example, if $\Omega$ consists of one binary operation:
\[
(a_1,a_2)\to a_1\circ a_2,\,a_1,a_2\in A\in{\mathfrak V},
\]
then we obtain the class of all ``magmas'' or ``groupoids'' (i.e. nonassociative semigroups).
The free object $F(X)$ consists of all words with parentheses $(x_{i_1}\cdots x_{i_k})\circ(x_{i_{k+1}}\cdots x_{i_n})$.

\begin{definition}\label{definition of Omega-identities}
{\rm (i)} Let $\mathfrak V$ be the class of all $\Omega$-algebras and let $F({\mathfrak V})$ be the free $\Omega$-algebra freely generated by the set $X=\{x_1,x_2,\ldots\}$.
For
\[
u(x_1,\ldots,x_n),v(x_1,\ldots,x_n)\in F({\mathfrak V})
\]
the $\Omega$-algebra $A\in{\mathfrak V}$ satisfies the identity
$u=v$ if
\[
u(a_1,\ldots,a_n)=v(a_1,\ldots,a_n)
\]
for all $a_1,\ldots,a_n\in A$. The class $\mathfrak W$ of all $\Omega$-algebras $A\in{\mathfrak V}$
satisfying a given system of identities $\{u_i=v_i \mid i\in I\}$ is called the variety defined by this system.

{\rm (ii)} Let $X$ be a countable set and let ${\rm Id}({\mathfrak W})$ be the set of all identities of $\mathfrak W$ in $F(X)$.
Then the factor-algebra $F({\mathfrak W})=F(X)/{\rm Id}({\mathfrak W})$
of $F(X)$ modulo the relations $u=v$ from ${\rm Id}({\mathfrak W})$ is the relatively free algebra of countable rank in $\mathfrak W$.
Similarly one defines the relatively free $\Omega$-algebra $F_d({\mathfrak W})$ of a finite rank $d$.
\end{definition}

Important examples of varieties are the classes of all groups, all semigroups, all commutative and associative algebras,
all associative, Lie, Jordan or all nonassociative algebras over a field $K$.
The corresponding free objects of countable rank are the free group $G(X)$, free semigroup $S(X)$, the polynomial algebra $K[X]$, the free associative algebra $K\langle X\rangle$,
the free Lie algebra $L(X)$, the free Jordan algebra $J(X)$, the free nonassociative algebra $K\{X\}$ and their ring counterparts
${\mathbb Z}[X]$, ${\mathbb Z}\langle X\rangle$, etc.

There are also other important classes of nonassociative algebras.

Alternative algebras are defined by the condition that every two elements generate an associative subalgebra, i.e. they satisfy the polynomial identities
\[
(x_1,x_1,x_2)=(x_1,x_2,x_1)=(x_2,x_1,x_1)=0,
\]
where $(x_1,x_2,x_3)=(x_1x_2)x_3-x_1(x_2x_3)$ is the associator.
The class of alternative algebras was implicitly discovered through the algebra of octonions in 1843
by John T. Graves and slightly later independently by Arthur Cayley.
But it was formally defined and isolated as a distinct algebraic variety by Max Zorn \cite{Zo} in 1930,
who credited key insights to Emil Artin.

(Right) Leibniz algebras satisfy the Leibniz identity
\[
[[x_1,x_2],x_3]=[x_1,[x_2,x_3]]+[[x_1,x_3],x_2],
\]
i.e. the right multiplication by any element is a derivation.
Leibniz algebras are a natural generalization of the Lie algebras. They were discovered in 1965 by Blokh \cite{Blo}, who called them D-algebras.
They attracted interest after the paper by Loday \cite{Lod}.

Left-commutative algebras satisfy the identity
\[
x_1(x_2x_3)=x_2(x_1x_3).
\]
Right-symmetric algebras are defined by the identity
\[
(x_1,x_2,x_3)=(x_1,x_3,x_2)
\]
and similarly one defines the right-commutative and left-symmetric algebras.
The first example of a left-commutative algebra and right-symmetric algebra is the Witt algebra in one variable
\[
W_1^{\rm{rsym}}=\left\{f\frac{d}{dx}\mid f\in K[x]\right\}
\]
with multiplication
\[
\left(f_1\frac{d}{dx}\right)\ast\left(f_2\frac{d}{dx}\right)
=\left(f_2\frac{df_1}{dx}\right)\frac{d}{dx}
\]
which appeared already in the paper by Cayley \cite{Ca} in 1857.
Cayley also found a realization of the algebras $W_d^{\rm{rsym}}$ in terms of rooted trees.

Right-symmetric and left-commutative algebras are called Novikov algebras.
The opposite algebras of Novikov algebras (i.e. left-symmetric and right-commutative algebras)
were introduced by Gel'fand and Dorfman \cite{GD} giving
an algebraic approach to the notion of Hamiltonian
operator in finite-dimensional mechanics and the formal calculus of
variations. Independently, later Novikov algebras were rediscovered by Balinskii and Novikov \cite{BN}
in the study of equations of hydrodynamics. The Witt algebra in one variable considered above is also the first example of a Novikov algebra.

Bicommutative algebras are nonassociative algebras which satisfy the polynomial identities of left- and right-commutativity
\[
x_1(x_2x_3)=x_2(x_1x_3)\;{\rm and}\;(x_1x_2)x_3=(x_1x_3)x_2.
\]

Finally, dual Leibniz algebras, called also Zinbiel algebras (Leibniz written in reverse order) satisfy the identity
\[
(x_1x_2)x_3 = x_1(x_3x_2 + x_2x_3).
\]

If the variety $\mathfrak W$ has the same identities as an algebra $A\in\mathfrak V$, then it is generated by $A$.
This definition agrees with the following theorem of Bikhoff.

\begin{theorem}
\label{HSB-theorem}{\rm (HSP-Theorem of Birkhoff \cite{Birk1})}
{\rm (i)} The class $\mathfrak W$ of algebraic systems from $\mathfrak V$ is a variety
if and only if it is closed under homomorphic objects ($\mathcal H$), subobjects ($\mathcal S$), and cartesian products ($\mathcal P$).

{\rm (ii)} The variety ${\mathfrak W}={\rm var}(A)$ is generated by $A$ if and only if
\[
{\mathfrak W}={\mathcal H}{\mathcal S}{\mathcal P}(A).
\]
\end{theorem}

\begin{definition}\label{definition of basis of identities}
{\rm (i)} Any set of defining identities of the variety ${\mathfrak W}$ of $\Omega$-algebras is called a basis of the identities of $\mathfrak W$.
The identity $u=v$ follows from the identities $\{u_i=v_i\mid i\in I\}$ if it holds for all $\Omega$-algebras which satisfy the identities
$u_i=v_i$, $i\in I$.

{\rm (ii)} The variety $\mathfrak W$ has a finite basis (or is finitely based) if it can be defined by a finite system of identities.
Otherwise it is infinitely (or non-finitely) based.
\end{definition}

The main topic of our paper is the following problem.
\begin{problem}\label{finite basis problem} {\rm (The finite basis (or Specht) problem)}
{\rm (i)} Let $\mathfrak V$ be the class of all algebraic system of a given signature $\Omega$.
Is every variety ${\mathfrak W}\subset\mathfrak V$ finitely based?

{\rm (ii) (Weaker version):} If the variety satisfies an explicitly given identity, does it have a finite basis of identities?

{\rm (iii)} If non-finitely based varieties exist, give an explicit example.

{\rm (iv)} What happens with the identities of interesting and important algebras?
\end{problem}

The origin of the problem is in the paper by Neumann \cite{Neu} from 1937. Translated in the modern language, the problem in \cite{Neu} is the following.

\begin{problem}\label{the problem of Bernard Neumann}
Is every variety generated by a finite group finitely based?
\end{problem}

In 1950 Wilhelm Specht \cite{Sp} asked the finite basis problem for associative algebras.
Again, the problem can be restated in the modern language:

\begin{problem}\label{original Specht problem}
Is it possible to define any variety of unitary associative algebras over a field of characteristic $0$
by a finite number of polynomial identities?
\end{problem}

Since the paper \cite{Sp} is in German, below we give some details.

Theorem 10 in \cite{Sp} states that for a unitary associative algebra $R$ over a field of characteristic 0
there is a chain of proper multilinear polynomial identities
\[
L_k(x_1,\ldots,x_k),L_{k+1}(x_1,\ldots,x_{k+1}),\ldots,
\]
such that all polynomial identities of $R$ follow from this chain.
(The multilinear polynomial $L_n(x_1,\ldots,x_n)$ is proper if $L_n(x_1,\ldots,x_{i-1},1,x_{i+1}\ldots,x_n)=0$ for all $i$.)

If the minimal degree of the polynomial identities of $R$ is equal to $k$, then Theorem 11 in \cite{Sp} gives that
there exists a reduced chain of proper multilinear polynomials
\[
L_k(x_1,\ldots,x_k),L_{k_1}(x_1,\ldots,x_{k_1}),\ldots,\,k<k_1<\cdots,
\]
such that all polynomial identities of $R$ follow from this chain and for any $n$
the $S_{k_n}$-module generated by $L_{k_n}(x_1,\ldots,x_{k_n})$ does not contain cconsequences of the identities of lower degree.

Here is the translation of the Specht problem as stated in \cite{Sp}:

\begin{problem}\label{original Specht problem-2}
In this case, I have not been able to decide whether there are rings with minimal degree of the polynomial identities equal to $k$
with an infinitely long reduced chain of proper multilinear identities, or whether every reduced chain of a ring terminates at a finite point.
\end{problem}

The influence of the paper by Specht \cite{Sp} and the activity trying to solve the problem
were so strong that now a variety of any algebraic systems is called Spechtian
if it and all its subvarieties can be defined by finite systems of identical relations.

In the language of varieties of algebraic systems a weaker version of Problem \ref{naive main problem} for classification
of all algebraic systems, naively stated again, is the following.

\begin{problem}\label{naive problem for varieties}
Describe all varieties of groups, semigroups, rings and algebras over a field (both associative and nonassociative).
\end{problem}

If the answer to the finite basis problem is negative then the idea to describe all varieties does not work.

In the rest of the paper we discuss various aspects of the Specht problem.
Since starting from the 1960s there is a lot of activity on the topic, the list of references in the paper is very far from being complete.

\section{Trivial answer: Infinitely based varieties of nonassociative objects}
The following theorem gives an easy example of varieties which are not finitely based.

\begin{theorem}\label{trivial eaxamples}
The varieties

$\bullet$ of binary magmas with zero (with $0$ and one nonassociative binary operation, i.e. nonassociative semigroups with zero);

$\bullet$ of nonassociative rings;

$\bullet$ of nonassociative algebras over any field $K$

\noindent defined by the identities
\[
u_n=(x_1x_2)(x_3x_4)x_5\cdots x_{n-2}(x_{n-1}x_n)=0,\quad n=6,7,\ldots,
\]
are not finitely based.
(The parentheses are left normed, e.g., $xyz=(xy)z$.)
\end{theorem}

The proof is almost obvious. For example, if we work in the free nonassociative algebra $K\{ X\}$ over a field, then all consequences of $u_n=0$
are linear combinations of monomials obtained by replacing some of the variables $x_1,\ldots,x_n$ by monomials and multiplication from the left and the right by monomials.
Clearly, in this way we cannot obtain a monomial $u_m$ for $m>n$.

\section{Meeting points with mathematical logic}

Obviously the finite basis problem can be restated in the language of mathematical logic:

\begin{problem}\label{translation in the language of logic}
Is the system of the identities of the variety finitely axiomatizable?
\end{problem}

The following problem was asked by Tarski \cite{Tar}.

\begin{problem}\label{problem of Tarski}
Is there an algorithm which, given an
arbitrary finite algebra $A$ of finite type (i.e. with a finite number of operations) as input, determines whether the variety ${\rm var}(A)$ is finitely based?
\end{problem}

Tarski attributed the problem to Perkins \cite{Per1}.
The negative answer to the problem was given by McKenzie \cite{McKe2}.
It uses the fact that there are Turing machines for which the halting problem is undecidable.

\begin{theorem}\label{Specht problem and Turing machine}
There is a construction which produces, for every Turing machine $\mathcal T$,
an algebra $A({\mathcal T})$ (finite and of finite type)
such that the Turing machine halts if and only if the algebra has a finite basis for its identities.
\end{theorem}

For more results related with logic see also Perkins \cite{Per2}, McNulty \cite{McNu} and McKenzie \cite{McKe1}.
To the best of our knowledge Problem \ref{problem of Tarski} is still open for finite semigroups.

\section{Positive results}

There are many methods which are used to prove that given varieties are finitely based.
In this section we shall survey several such methods.

\subsection{Structure theory of groups and algebras}\label{Specht property for nice finite objects}

A typical result using structure theory of finite objects is the following theorem.

\begin{theorem}\label{Finite objects}
Let $\mathfrak V$ be a variety generated by a finite group, associative, Lie, Jordan or alternative ring
or by a finite-dimensional associative, Lie, Jordan or alternative algebra over a finite field.
Then $\mathfrak V$ is finitely based.
\end{theorem}

This result was obtained by:

$\bullet$ Sheila Oates and Powell \cite{OaPo} for groups;

$\bullet$ Kruse \cite{Kru} and L'vov \cite{Lv1} for finite associative rings and algebras;

$\bullet$ Bahturin and Olshanski\v{\i} \cite{BaOl1} for finite Lie rings and algebras;

$\bullet$ Medvedev \cite{Med1} for finite Jordan rings and algebras;

$\bullet$ L'vov \cite{Lv3} for finite alternative rings and algebras.

The proofs of Theorem \ref{Finite objects} use the same idea which comes from the paper by Oates and Powell \cite{OaPo}.
(For Lie rings the proof involves also quasi-identities.)

\begin{definition}\label{critical groups and Cross varieties}
A finite group is called critical if it does not belong to the variety generated by its proper subgroups and factor groups.

A variety of groups $\mathfrak V$ is Cross if:

$\bullet$ $\mathfrak V$ has a finite basis of its identities.

$\bullet$ Finitely generated groups in $\mathfrak V$ are finite.

$\bullet$ $\mathfrak V$ contains a finite number of critical groups.

\noindent The definitions for rings and algebras are similar.
\end{definition}

Then Theorem \ref{Finite objects} is a consequence of the following result.

\begin{theorem}\label{finite generators and Cross}
A variety of groups, associative, Lie, Jordan or alternative rings and algebras is Cross if and only if it is generated by a finite object.
\end{theorem}

\subsection{The Higman-Cohen method}

The method is based on considerations in the spirit of the following result
which is a consequence of the result of Higman \cite{Hig1}:

\begin{theorem}\label{theorem of Higman}
Consider the finite sequences of nonnegative integers with partial ordering
\[
(a_1,\ldots,a_m)\preceq(b_1,\ldots,b_n)
\]
if there is a subsequence $(b_{i_1},\ldots,b_{i_m})$ of $(b_1,\ldots,b_n)$ such that
\[
a_1\leq b_{i_1},\ldots,a_m\leq b_{i_m}.
\]
Then the set of such sequences satisfies the descending chain condition and any set of sequences has a finite subset of minimal elements.
\end{theorem}

For the first time this was applied by Cohen \cite{Coh} who proved:

\begin{theorem}\label{theorem of Cohen}
The metabelian variety of groups satisfies the Specht property.
\end{theorem}

Most of the positive results on the Specht problem in 1970s and 1980s were obtained using this method.
Here are some typical results:

$\bullet$ Group theory: Vaughan-Lee \cite{V-L1}, Susan McKay \cite{McKa1}, \cite{McKa2}, Bryant and Newman \cite{BrNew}, Galina Sheina \cite{Shei1}.

$\bullet$ Lie algebras: Vaughan-Lee \cite{V-L4}, Bryant and Vaughan-Lee \cite{BrV-L}, Sheina \cite{Shei2};

$\bullet$ Associative algebras (the main tool for establishing positive results before the work of Kemer):
Latyshev \cite{La1}, \cite{La2}, Genov \cite{Gen1}, Popov \cite{Pop};

$\bullet$ Alternative algebras: Umirbaev \cite{Um};

$\bullet$ Monoids, i.e. semigroups with 1: Lee \cite{Lee1};

$\bullet$ Group graded Lie algebras: Daniela Martinez Correa and Koshlukov \cite{MCK}, Martinez Correa and Yasumura \cite{MCY}.

Vaughan-Lee \cite{V-L1} showed that the variety ${\mathfrak A}{\mathfrak N}_c\cap{\mathfrak N}_c{\mathfrak A}$
of abelian-by-(nilpotent of class $c$) and (nilpotent of class $c$)-by-abelian groups defined by the identities
\[
((x_1,\ldots,x_{c+1}),(x_{c+2},\ldots,x_{2c+2}))=((x_1,x_2),\ldots,(x_{2c+1},x_{2c+2}))=1
\]
satisfies the Specht property. Here $(u,v)=u^{-1}v^{-1}uv$ is the group commutator and the longer commutators are left normed: $(u,v,w)=((u,v),w)$.
McKay \cite{McKa2} established the Specht property for the center-by-metabelian variety of groups defined by the identity
$(((x_1,x_2),(x_3,x_4)),x_5)=1$.
Bryant and Newman \cite{BrNew} generalized this for the variety ${\mathfrak N}_c{\mathfrak A}\cap{\mathfrak N}_2{\mathfrak N}_c$.
In particular, this implies that the variety ${\mathfrak N}_2{\mathfrak A}$ is also Spechtian.
All these results were further generalized by Sheina \cite{Shei1}.

Vaughan-Lee \cite{V-L4} proved that over a field of characteristic different from 2 the variety of center-by-metabelian Lie algebras satisfies the Specht property.
(Later we shall see that over a field of characteristic 2 Vaughan-Lee gave an example of a not finitely based variety of center-by-metabelian Lie algebras.
This was generalized by Bryant and Vaughan-Lee \cite{BrV-L} for the variety of Lie algebras ${\mathfrak N}_2{\mathfrak A}$.
Sheina \cite{Shei2} went further and handled the case of the variety ${\mathfrak N}_c{\mathfrak A}\cap{\mathfrak N}_2{\mathfrak N}_c$. (Again ${\rm char}(K)\not=2$.)

Latyshev \cite{La1} and Genov \cite{Gen1}, \cite{Gen2} proved that over a field of characteristic 0 the variety of associative algebras ${\mathfrak N}_c{\mathfrak A}$
defined by the polynomial identity
\[
[x_1,x_2]\cdots[x_{2c-1},x_{2c}]=0
\]
is Spechtian. (By the theorem of Yu.N. Maltsev \cite{Malts}
this variety is generated by the algebra of $c\times c$ upper triangular matrices.)
Latyshev  \cite{La2} and Popov \cite{Pop} established the finite basis property for the variety of associative algebras over a field of characteristic 0
defined by the identity
\[
[x_1,x_2,x_3]\cdots[x_{3c-2},x_{3c-1},x_{3c}]=0.
\]

Umirbaev \cite{Um} proved that over a field of characteristic different from 2 and 3 the variety of solvable alternative algebras satisfies the Specht property.
The proof uses a result of Pchelintsev \cite{Pch1} that in characteristic different from 2 and 3 any solvable variety of alternative algebras is a subvariety of a variety
${\mathfrak N}_c{\mathfrak A}\cap{\mathfrak N}_3{\mathfrak N}_c$, i.e. for a solvable alternative algebra $R$ one has $(R^2)^c=(R_c)^3=0$. This is an analogue of the result for Lie algebras due to
Bryant and Vaughan-Lee \cite{BrV-L} and Sheina \cite{Shei2}.

Let $R$ be an algebra and let $G$ be a group. The algebra $R$ is $G$-graded if there is a vector space decomposition as a direct sum
\[
R=\bigoplus_{g\in G}R_g\;{\rm and}\;R_gR_h\subseteq R_{gh},\;g,h\in G.
\]
In the special case when $R=M_c(K)$ is the $c\times c$ matrix algebra there is a canonical ${\mathbb Z}_c$-grading:
the matrix units $e_{ij}$ are of degree $j-i\;{\rm (mod}\;c)$. If $G$ is an abelian group, then the elementary $G$-grading of the algebra $UT_c(K)$
of $c\times c$ upper triangular matrices is induced by the assigning of degree $g_i\in G$ to the matrix units $e_{i,i+1}$, $i=1,2,\ldots,c-1$.
By analogy with the ordinary case one studies the $G$-graded polynomial identities of the $G$-graded algebra $R$.
Martinez Correa and Koshlukov \cite{MCK} showed that if ${\rm char}(K)=0$ or ${\rm char}(K)=p\geq c$, then the variety of ${\mathbb Z}_c$-graded Lie algebras
generated by $UT_c^{(-)}(K)$ with the canonical grading satisfies the Specht property. A similar result was established by Martinez Correa and Yasumura \cite{MCY} for the Lie algebra $UT_3^{(-)}(K)$,
${\rm char}(K)\not=2$, equipped with an elementary grading by an arbitrary abelian group $G$.

\subsection{Methods of commutative algebra}
It is well known, see e.g. \cite[Example 4.35, page 132]{Ba3} that over an infinite field $K$ the polynomial identities
of the Lie algebra $UT_c^{(-)}(K)$ of $c\times c$ upper triangular matrices follow from the identity
\[
[[x_1,x_2],\ldots,[x_{2c-1},x_{2c}]]=0,
\]
i.e. ${\mathfrak N}_{c-1}{\mathfrak A}={\rm var}(UT_c^{(-)}(K))$.
Over a field of characteristic 0 Krasilnikov \cite{Kra} proved the Specht property for ${\rm var}(UT_c^{(-)}(K))$.

\begin{theorem}\label{Krasilnikov, Lie upper triangular}
Let ${\rm char}(K)=0$. Then the variety of Lie algebras ${\mathfrak N}_{c-1}{\mathfrak A}$ satisfies the finite basis property.
\end{theorem}

The proof is based on a combination of several methods. Applying representation theory of the symmetric and general linear groups
in the spirit of the approach developed by the author \cite{Dr3} and Berele \cite{Be}
and the fact that $\dim(UT_c^{(-)}(K))<\infty$, one derives that modulo the identity
\[
[[x_1,x_2],\ldots,[x_{2c-1},x_{2c}]]=0
\]
every polynomial identity is equivalent to an identity of $\leq \dim(UT_c^{(-)}(K))$ variables.
Then the finite basis problem is reduced to the problem for the finite generation of a finitely generated module of the tensor product
\[
K[Y_n]^{S_n}\otimes\cdots\otimes K[Y_n]^{S_n}
\]
of $c$ copies of the algebra of symmetric polynomials $K[Y_n]^{S_n}$ of $n=\dim(UT_c^{(-)}(K))$ variables.
Clearly, now the proof follows from the Hilbert Basis Theorem. Modification of the proof of Krasilnikov works also for associative algebras.
In this way one can obtain an alternative proof of the results of Latyshev \cite{La1} and Genov \cite{Gen1}, \cite{Gen2} for the Specht property of the variety
${\rm var}(UT_c(K))$ over a field $K$ of characteristic 0.

Hilbert Basis Theorem was used also by McKay \cite{McKa2}. She showed that if $\mathfrak V$ is a subvariety of the variety ${\mathfrak N}_2{\mathfrak A}$
of groups defined by a system of identities in $n$ variables, then $\mathfrak V$ is finitely based.

Let $\mathfrak B$ be the variety of all bicommutative algebras.
Dzhumadil'daev, Ismailov, and Tulenbaev \cite{DIT} (see also \cite{DT1})
described the free bicommutative algebra $F({\mathfrak B})$ of countable rank and its main numerical invariants.
By \cite{DT1} the square $F^2({\mathfrak B})$ of the algebra $F({\mathfrak B})$ is a commutative and associative algebra.
In our paper with Zhakhayev \cite{DrZh} we used that
\[
F^2({\mathfrak B})\cong\omega(K[Y])\omega(K[Z])\subset K[Y,Z],
\]
where $\omega(K[Y])$ and $\omega(K[Z])$ are the augmentation ideals of the corresponding polynomial algebras. The following theorem was one of the main results of \cite{DrZh}.

\begin{theorem}\label{Specht for bicommutative algebras}
Over an arbitrary field $K$ of any characteristic the variety of bicommutative algebras satisfies the Specht property.
\end{theorem}

In the case of characteristic 0 the proof uses representation theory of the symmetric and general linear group to show that
modulo the left- and right-commutativity the polynomial identities are equivalent to identities in two variables.
Then the result follows immediately from the Hilbert Basis Theorem. In the case of fields of positive characteristic we applied the Higman-Cohen method.

In our paper with Ismailov, Mustafa and Zhakhayev \cite{DrIMZh} we introduced bicommutative superalgebras. With methods similar to those in \cite{DrZh}
we obtain similar results for the Specht property of varieties of bicommutative superalgebras.

\subsection{The role of nilpotency}

One of the key results in the theory of associative algebras is the Dubnov-Ivanov-Nagata-Higman theorem \cite{DuI,Nag,Hig2}.

\begin{theorem}
Let the associative algebra $R$ is nil of bounded index $n$, i.e. $R$ satisfies the identity $x^n=0$.
If the field $K$ is of characteristic $0$ or $p>n$, then $R$ is nilpotent of class $N$ for some $N$ depending on $n$.
\end{theorem}

Dzhumadil'daev and Tulenbaev \cite{DT2} proved an analogue of this theorem for Zinbiel (or dual Leibniz) algebras.

\begin{theorem}
Over a field of arbitrary characteristic every Zinbiel nil-algebra of nil-index $n$ is nilpotent of class $\leq 2^n-1$.
\end{theorem}

The authors of \cite{DT2} established also other results for solvability, nil and nilpotency of Zinbiel algebras. In the proofs they used essentially that Zinbiel algebras are right-commutative.

An important moment in the proofs of most of the results in this subsection is that one has detailed information for the vector space bases of the relatively free algebras in the varieties.
This was also essential for the results of Naurazbekova and Umirbaev \cite{NaUm} on the variety of Zinbiel algebras.

\begin{theorem}
Let $K$ be a field of characteristic $0$.

{\rm (i)} The variety $\mathfrak Z$ of Zinbiel algebras is generated by its relatively free algebra $F_1({\mathfrak Z})$ of rank $1$.

{\rm (ii)} Any proper subvariety of $\mathfrak Z$ is nilpotent.

{\rm (iii)} The variety $\mathfrak Z$ satisfies the Specht property.
\end{theorem}

As a consequence of these results, in the case of characteristic 0, \cite{NaUm} gives also a new proof of the result of \cite{DT2}
that finite-dimensional Zinbiel algebras are nilpotent.

For an algebra $R$ one defines the left associator $(R,R,R)_n$ inductively:
\[
(R,R,R)_1=(R,R,R)=\{(r_1,r_2,r_3)\mid r_1,r_2,r_3\in R\},
\]
\[
(R,R,R)_n=(R,R,(R,R,R)_{n-1}),\,n\geq 2.
\]
Dotsenko, Ismailov and Umirbaev, \cite{DoIsUm} established a special kind of nilpotency for Novikov algebras.
(The authors of \cite{DoIsUm} consider the opposite algebras of Novikov algebras, as introduced by Gel'fand and Dorfman \cite{GD}.
We shall translate their results in the language of Novikov algebras.)

\begin{theorem}
Over a field of characteristic zero a Novikov algebra $R$ which satisfies a nontrivial identity
is left associator nilpotent, i.e. $(R,R,R)_n=0$ for some $n$.
\end{theorem}

The main result of \cite{DoIsUm} is the following:

\begin{theorem}
Over a field of characteristic zero the variety of Novikov algebras satisfies the Specht proeprty.
\end{theorem}

The proofs of most of the results in \cite{DoIsUm} are in the language of operads.
The paper gives also a comprehensive survey of results on identities of Novikov algebras
and contains other interesting results, especially concerning the Specht property of commutative associative differential algebras.
We shall state only one more result from \cite{DoIsUm}.
Right-symmetric algebras algebras are also called pre-Lie algebras. They become Lie algebras with respect to the commutator $[x,y]=xy-yx$.

\begin{corollary}
If a Novikov algebra over a field of characteristic zero satisfies a nontrivial identity, then it is Lie solvable.
\end{corollary}

\subsection{Weak polynomial identities}

Let $R$ be an associative algebra generated by a vector subspace $V$. Then the polynomial
$f(x_1,\ldots,x_n)\in K\langle X\rangle$ is a weak polynomial identity for the pair $(R,V)$ if
\[
f(v_1,\ldots,v_n)=0\;{\rm for}\;{\rm all}\;v_1,\ldots,v_n\in V.
\]
The consequences of the weak identities $f_i(x_1,\ldots,x_{n_i})$, $i\in I$, are the elements of the ideal generated by all
$f_i(u_1,\ldots,u_{n_i})$ where $u_1,\ldots,u_{n_i}$ are linear combinations of the elements of $X$.
If $V$ is a Lie subalgebra of $R$, then for the consequences it is allowed to replace $x_1,\ldots,x_{n_i}$ by the elements of the Lie algebra generated by $X$
in $K\langle X\rangle$. Weak polynomial identities were introduced by Razmyslov \cite{Raz1}, \cite{Raz2}, who applied them to construct central polynomials for $M_c(K)$
and to find the basis of the polynomial identities in $M_2(K)$ and $sl_2(K)$, ${\rm char}(K)=0$.

The following theorem was proved by Volichenko \cite{Vol2}.

\begin{theorem}\label{Volichenko week PI}
Over a field of characteristic $0$ the variety of Lie algebras
${\mathfrak A}{\mathfrak N}_2$ defined by the identity $[[x_1,x_2,x_3],[x_4,x_5,x_6]]=0$ is Spechtian.
\end{theorem}

The key moment of the proof is the Specht property for the weak polynomial identity $[x,y,z]=0$.


\subsection{Structure theory of T-ideals -- the approach of Kemer}

In the case of characteristic 0 Kemer \cite{Kem1}, \cite{Kem2}, \cite{Kem3} (see also the books by Kanel-Belov and Rowen \cite{K-BR}
and Kanel-Belov, Karasik and Rowen \cite{K-BKR} for the exposition)
developed structure theory of T-ideals of the free associative algebra
in the spirit of ideal theory in polynomial algebras, combined with applications of the theory of superalgebras.
One of the important results in the theory of Kemer is the positive solution of the Specht problem.

\begin{theorem}\label{Kemer for associative algebras}
Over a field of characteristic $0$ every variety of associative algebras has a finite basis of its polynomial identities.
\end{theorem}

Analogs of the results of Kemer were obtained for $G$-graded associative algebras and associative algebras with involution.
The following theorems were established independently by Irina Sviridova \cite{Sv1} for $G$-graded algebras for finite abelian groups and
by Aljadeff and Kanel-Belov \cite{AK-B} for any finite groups and by Sviridova \cite{Sv2} for algebras with involution.

\begin{theorem}
Let $G$ be a finite group and let $K$ be a field of characteristic $0$. Then every variety of $G$-graded associative algebras over $K$ is finitely based.
\end{theorem}

\begin{theorem}
Let $K$ be a field of characteristic $0$. Then every variety of associative algebras with involution over $K$ is finitely based.
\end{theorem}

The Specht problem is still open for varieties of Lie algebras over a field of characteristic 0.

\section{Finite algebraic system with a finite system of operations}

In Subsection \ref{Specht property for nice finite objects} we commented that the Specht problem has a positive solution for nice finite objects --
groups and associative, Lie and Jordan rings and algebras. In 1948 Birkhoff \cite{Birk2} asked the following question.

\begin{question}\label{question of Birkhoff-1948}
Does every finite algebraic system with a finite system of operations possess a
finite set of identities from which all others are derivable?
\end{question}

The answer into affirmative for systems with two elements was given by Lyndon \cite{Ly1} in 1951.

\begin{theorem}\label{systems with two elements}
Every $\Omega$-algebra with two elements and a finite set of operations has a finite basis of its identities.
\end{theorem}

As Lyndon commented ``it is perhaps surprising
that this problem is not entirely trivial even for algebras that contain only two elements.''

The first counterexample to Question \ref{systems with two elements} was found by Lyndon \cite{Ly2} in 1954.

\begin{theorem}\label{counterexample of Lyndon}
There is an algebraic system which consists of $7$ elements and has one binary operation which does not possess
a finite system of identities.
\end{theorem}

The magma $M=\{0,e,b_1,b_2,c,d_1,d_2\}$ of Lyndon has 0 and its only
nonzero products are
\[
ce=c,cb_j=d_j,d_je=d_j,d_jb_k=d_j,\quad j,k=1,2.
\]
Its identities, with left normed notation, e.g. $xyz=(xy)z$, follow from:
\[
0x=0,x0=0,x(yz)=0,
\]
\[
x_1x_2\cdots x_nx_1=0,\quad n=1,2,\ldots,
\]
\[
x_1x_2\cdots x_nx_2=x_1x_2\cdots x_n,\quad n=1,2,\ldots.
\]

The variety ${\rm var}(M)$ generated by $M$ was studied by Lee \cite{Lee5}.

\begin{theorem}\label{paper Lee5}
{\rm (i)} The subvarieties of ${\rm var}(M)$ form a chain with five elements, all of which
except ${\rm var}(M)$ are Cross varieties.

{\rm (ii)} The magma $M$ has a quotient of $6$ elements which also generates $\mathfrak M$.
\end{theorem}

The result of Lyndon naturally raises a question.

\begin{problem}\label{minimal counterexample for finite objects}
What is the minimal number of elements of a finite algebraic system $M$ without a finite basis of identities?
\end{problem}

A result of Vishin reduced Problem \ref{minimal counterexample for finite objects} to {\it what is the minimal number of elements -- $3$ or $4$?}

\begin{theorem}\label{4 elements} {\rm (Vishin \cite{Vi})}
Let $M=\{0,1,2,3\}$ be the four-element binary magma with the only nonzero products
\[
1\ast 2=1,3\ast 3=3,3\ast 2=1.
\]
Then $M$ is not finitely based.
\end{theorem}

The magma of Vishin satisfies the identities
\[
x0=0x=x_1(x_2(x_3x_4))=0
\]
and the left normed identities
\[
x_1x_2x_2=x_1x_2,
\]
\[
x_1x_2x_3x_4\cdots x_nx_2=x_1x_3x_2x_4\cdots x_nx_3,n\geq 3,
\]
\[
x_1\cdots x_kx_{k+1}\cdots x_nx_2=x_1\cdots x_{k+1}x_k\cdots x_nx_2,k=4,\ldots,n-1,n\geq 5,
\]
\[
x_1x_2\cdots x_nx_1=0,n\geq 2.
\]

Finally, the answer to Problem \ref{minimal counterexample for finite objects} was given in 1965 by Murskii \cite{Mur1}.
He constructed the following magma $M$:
\begin{center}
\includegraphics[width=2.1cm]{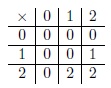}
\end{center}

\begin{theorem}\label{Murski n=3}
The three-element binary magma $M$ of Murskii is not finitely based.
\end{theorem}

The proof follows from the fact that
for $n\geq 3$ the magma $M$ satisfies the identity
\[
x_1(x_2(x_3\cdots (x_{n-1}(x_nx_1))\cdots))=(x_1x_2)(x_n(x_{n-1}\cdots (x_4(x_3x_2))\cdots ))
\]
which does not follow from any system of identities in less than $n$ variables.

Je\u{z}ek \cite{Je} found three more non-finitely based three-element magmas. All they are idempotent, i.e. satisfy the identity $x^2=x$
and two of them are commutative.

In a series of papers Murskii \cite{Mur2}, \cite{Mur3}, \cite{Mur4} studies the finite basis problem for finite algebraic systems
with methods typical for mathematical logic. The papers are devoted to the following problem.

\begin{problem}\label{how many are counterexamples}
{\rm (i)} For a fixed finite signature $\Omega$ how many are the $k$-element $\Omega$-algebras without a finite basis of their identities?

{\rm (ii)} For a fixed number $k$ of elements of the algebra how big is the number of $\Omega$-algebras without a finite basis of their identities when
the number and the arity of the operations in $\Omega$ increases?
\end{problem}

Let $\Omega$ be a finite set of operations which contains also nonunary operations (i.e. $n$-ary operations for $n>1$).
Let $N(\Omega,k)$ be the set of all $\Omega$-algebras with $k$ elements and let $P(\Omega,k)$ be the subset of the $\Omega$-algebras without the finite basis property.
Let
\[
\delta(\Omega,k)=\frac{\vert P(\Omega,k)\vert}{\vert N(\Omega,k)\vert}
\]
be the ratio between the cardinality of $P(\Omega,k)$ and $N(\Omega,k)$. The results of Murskii \cite{Mur2} show that ``almost all'' finite algebras have a finite basis for their identities.
In \cite{Mur3} he makes it more precise for the case when $\Omega$ consists of one binary operation.

\begin{theorem}\label{theorems of Murski} {\rm(\cite{Mur2}, \cite{Mur3}) (i)} For a fixed $\Omega$
\[
\lim_{k\to\infty}\delta(\Omega,k)=0.
\]
For a fixed $k$ the ratio $\delta(\Omega,k)$ also tends to $0$ when $\vert\Omega\vert\to\infty$.

{\rm (ii)} When $\Omega$ consists of one binary operation, then $\delta(\Omega,k)={\mathcal O}(k^{-6})$.

{\rm (iii)} If $\Omega$ consists of one binary operation, then any finite $\Omega$-algebra which contains as a subalgebra the magma $M$ of Murskii from Theorem \ref{Murski n=3}
does not have a finite basis of its identities.
\end{theorem}

In \cite{Mur4} Murskii showed that the situation is different if no distinction is made between algebras determining
the same clone of operations on the common set of their elements, see \cite{Mur4} for details.

\section{Infinitely based varieties of semigroups}

\subsection{The first examples}
The first examples of varieties of semigroups which are not finitely based
were given by Birjukov \cite{Birj} in 1965 and Austin \cite{Au} in 1966.

\begin{theorem}\label{fist example semigroups}
The following systems of identities define varieties of semigroups which are not finitely based:

{\rm (i) (Birjukov \cite{Birj})}
\[
y(x_1\cdots x_n)z(yx_n\cdots x_1)yz=y(x_1\cdots x_n)z(x_n\cdots x_1)yz,n=1,2,\ldots
\]

{\rm (ii) (Austin \cite{Au})}
\[
((x_1\cdots x_n)y)^2=y(x_1\cdots x_n)y(x_n\cdots x_1),n=2,3,\ldots
\]
\end{theorem}

For more information concerning results on identities in semigroups before 1985 see the survey by Shevrin and Volkov \cite{SheV}.
(More than 40 examples of not finitely based varieties of semigroups can be found there.)
For more recent results see the survey by Volkov \cite{Volk2}.

\subsection{Finite semigroups with infinite bases of identities}
In 1969  Perkins \cite{Per3} constructed the first example of a finite semigroup without a finite basis of its identities:

\begin{theorem}\label{6-element semigroup of matrices}
The following six-element semigroup of $2\times 2$ matrices (known as the six-element Brandt monoid $B_2^1$) is not finitely based:
\[
\left\{ \left(\begin{array}{cc}
0 & 0 \\
0 & 0 \\
\end{array}\right),\left(\begin{array}{cc}
1&0\\
0&0\\
\end{array}\right),\left(\begin{array}{cc}
0&1\\
0&0\\
\end{array}\right),\left(\begin{array}{cc}
0&0\\
1&0\\
\end{array}\right),\left(\begin{array}{cc}
0&0\\
0&1\\
\end{array}\right),\left(\begin{array}{cc}
1&0\\
0&1\\
\end{array}\right)\right\}.
\]
\end{theorem}

Volkov, Goldberg and Kublanovsky \cite{VGK} found another six-element semigroup that does not have a finite basis of its identities, but nevertheless
generates a variety whose finite membership problem admits a polynomial algorithm.

Up to isomorphism or anti-isomorphism, there exist 15973 distinct semigroups with 6 elements
(and 1373 of them are monoids).
Four of them (2 of them are monoids) are not finitely based.
Each of the remaining 15969 semigroups with 6 elements has finitely based identities
(Lee and Li \cite{LeeLi} for monoids; Lee and Zhang \cite{LeeZh} for the other semigroups).

In his paper \cite{Tar} Tarski raised explicitly the finite basis problem for semigroups of order five or less
which attracted the interest of several mathematicians
(Bol'bot, Edmunds, Karnofsky, Tishchenko, Trahtman).
A solution to this problem was eventually completed by Trahtman \cite{Tra3}, see also Lee \cite{Lee6}:

\begin{theorem} Every semigroup with $\leq 5$ elements has a finite basis of its identities.
\end{theorem}

If a variety of algebraic systems satisfies an infinite system of independent identities
(i.e. removing one of the identities we obtain a larger variety)
then it has a continuum of subvarieties. (For the proof, apply the Cantor diagonal argument.)

\begin{theorem}\label{countably many varieties if semigroups}
{\rm (Trahtman \cite{Tra2}, Edmunds, Lee and Lee \cite{ELL})}
There exist six-element semigroups which generate a variety with a continuum of subvarieties.
\end{theorem}

The semigroup of Trahtman has an infinite basis of identities and that of Edmunds, Lee and Lee
is finitely based.

Up to isomorphism and anti-isomorphism, every semigroup with $\leq 5$ elements, with one possible exception,
generates a variety with at most countably many subvarieties. This possible exception is the semigroup given by the following multiplication table

\begin{center}
\includegraphics[width=3.2cm]{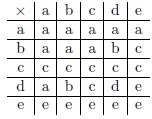}
\end{center}
\noindent
(for the proof see the paper by Lee \cite{Lee6}).

On the other hand, Trahtman \cite{Tra1} constructed a variety of semigroups with 1 without an infinite independent system of identities.

\begin{theorem}\label{no independent system for semigroups}
The variety of semigroups with $1$ defined by the identities
\[
yx^4y=yx^5y,\; yxz_1^2\cdots z_n^2x^6y=y^2xz_1^2\cdots z_n^2x^6y,\;n=1,2,\ldots,
\]
does not have an independent system of identities.
\end{theorem}

We continue with two recent exotic examples where the Specht property depends on the existence of an involution.

\begin{theorem}\label{with and without involution}
The following five-element semigroups given by their defining relations
generate Spechtian varieties of semigroups but non-Specht varieties of involution semigroups.

{\rm (i) (Gao, Zhang and Luo \cite{GZL})}:
\[
A_0^1=\langle e,f\mid e^2=e,f^2=f,ef=0\rangle,
\]
\[
e^{\ast}=f,f^{\ast}=e.
\]

{\rm (ii) (Lee \cite{Lee3})}:
\[
A_2=\langle a,e\mid a^2=0,aea=a,e^2=eae=e\rangle,
\]
\[
a^{\ast}=a,e^{\ast}=e.
\]
\begin{center}
\includegraphics[width=7.2cm]{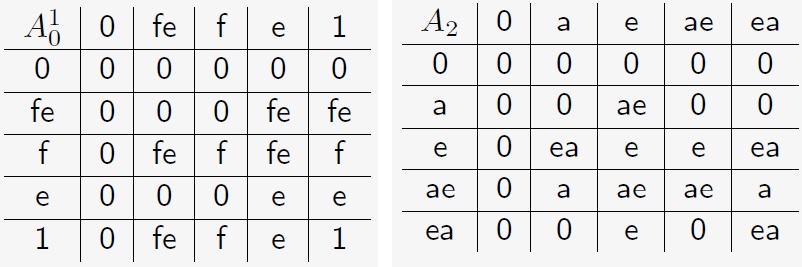}
\end{center}
\end{theorem}

The above result was further extended:

\begin{theorem}
Let $A_0$ be the four-element subsemigroup of the monoid $A_0^1$ in the previous theorem. Then

{\rm (i) (Lee \cite{Lee7})} As a semigroup $A_0$ generates a Spechtian variety.

{\rm (ii) (Lee \cite{Lee8})} As an involution semigroup $A_0$ generates a non-Spechtian variety
because contains a non-finitely based subvariety.

{\rm (iii) (Gao, Lee, Luo, Zhang \cite{GLLZ})}: As an involution semigroup $A_0$ has a finite basis of its identities.
\end{theorem}

Zhang, Gao, and Luo \cite{ZhGL} studied the Boolean matrix monoid $BR_n$ which consists of all $n\times n$ Boolean matrices with 1s on the main diagonal.
This monoid has two kinds of involutions -- the transposition and the skew transposition.
The authors solved completely the finite basis problem of $BR_n$ as an involution monoid with respect to the two involutions.
They also obtained other results concerning the identities of these involution monoids.

For a survey on other results concerning varieties of involution semigroups see the paper \cite{CrDol}
by Crvenkovi\'c and Dolinka and the book by Lee \cite{Lee4}. See also the paper by Jackson and Lee \cite{JL} for varieties with extreme properties
generated by finite monoids.

We conclude this section with results which connect the Specht problem for semigroups with graph theory and decision problems.
At the conference on Structural Theory of Automata, Semigroups
and Universal Algebra at the Universit\'e de Montr\'eal in 2003 Volkov introduced the
following problem.

\begin{problem}\label{NP-complete}
Does there exist a finite monoid $M$ such that the problem to determine
whether a finite monoid $M'$ belongs to the variety ${\rm var}(M)$ (the finite membership problem
for ${\rm var}(M)$) is NP-complete?
\end{problem}

This problem was inspired by a similar problem for finite general algebras solved by Sz\'ekely \cite{Sz} who constructed
a seven-element algebra $A$ such that the finite membership problem for ${\rm var}(A)$ is NP-complete.
Jackson and McKenzie \cite{JMcK} constructed a 55-element semigroup $S$ and 56-element monoid $S^1$ such that the finite membership problem
for the variety of semigroups ${\rm var}(S)$ and the variety of monoids ${\rm var}(S^1)$ are at least as difficult as determining if a finite
graph is 3-colorable. This implies that these membership problems are NP-hard. As a consequence the authors obtain that the variety ${\rm var}(S^1)$
has continually many subvarieties of monoids. They also give new interesting series of non-finitely based finite semigroups.

Quoting from \cite{JZ} ``This Brandt monoid is perhaps the most ubiquitous harbinger of complex behaviour in all finite semigroups.''
Recently Jackson \cite{J} proved that the variety of monoids generated by the six-element Brandt monoid $B_2^1$ has NP-hard variety membership problem.
By \cite{JL} the variety of monoids generated by $B_2^1$ does not have a finite basis of its identities
which is a monoid analogue of the semigroup results in Theorem \ref{6-element semigroup of matrices}.
Volkov \cite{Volk3} used the Brandt monoid to give a simple proof for the existence of a six-element semiring without a finite basis of its identities.
Gusev \cite{Gus3} showed that the variety of monoides generated by the 6-element Brandt monoid is finitely universal,
i.e. its lattice of subvarieties contains an isomorphic copy of every finite lattice.

\section{Infinitely based varieties of groups}

As we stated in Theorem \ref{Finite objects} the positive answer to Problem \ref{the problem of Bernard Neumann} as stated originally for finite groups
was given by Oates and Powell \cite{OaPo} in 1964. The first result which gives the existence of non-finitely based varieties of groups was obtained by Ol'shanskij \cite{Ol} in 1970.

\begin{theorem}\label{theorem of Olshanski}
There are continually many varieties of groups.
\end{theorem}

Since the free group $G(X)$, $X=\{x_1,x_2,\ldots\}$ has countably many elements, if all varieties of groups were finitely based,
then there would be only countably many varieties. Hence Theorem \ref{theorem of Olshanski}
immediately implies that there are infinitely based varieties. These varieties can be chosen to be solvable of class $5$
(i.e. are subvarieties of the variety ${\mathfrak A}^5$, where $\mathfrak A$ is the variety of abelian groups)
and of exponent $8pq$, where $p$ and $q$ are any coprime odd integers.
The proof of the theorem is based on the following lemma.

\begin{lemma}
Let $T_i$, $i=1,2,\ldots$, be a set of finite groups such that:

{\rm (i)} All groups $T_i$ belong to a locally finite variety $\mathfrak U$
(Locally finite = every finitely generated group is finite);

{\rm (ii)} All groups $T_i$ are monolithic (i.e. have a unique minimal normal subgroup which is contained in all proper normal subgroups);

{\rm (iii)} The group $T_i$ is not isomorphic to a factor of $T_j$, $j\not= i$.

If the prime number $p$ does not divide $\exp({\mathfrak U})$, then the product ${\mathfrak A}_p{\mathfrak U}$ has a continuum different subvarieties.
\end{lemma}

The first explicit examples of varieties of groups which were not finitely based were given in the same 1970 by Adyan \cite{Ad} and Vaughan-Lee \cite{V-L2}.

\begin{theorem}\label{example of Adyan} {\rm (Adyan \cite{Ad})}
For every odd $n\geq 4381$ the group identities
\[
u_p=(x^{pn},y^{pn})^n=1,\quad p\;prime,
\]
are independent (i.e., $u_p$ does not follow from the other identities).
Here $(u,v)=u^{-1}v^{-1}uv$ is the group commutator.
\end{theorem}

The proof uses methods from the solution of the Burnside problem by Novikov and Adyan \cite{NovAd}
which can be stated in the language of varieties:
{\it Is the variety of groups defined by the identity $x^k=1$ locally finite?}

\begin{theorem}\label{group example of Vaughan-Lee} {\rm (Vaughan-Lee \cite{V-L2})}
Let
\[
w_k=((x,y,z),(x_1,x_2),(x_3,x_4),\ldots,(x_{2k-1},x_{2k}),(x,y,z)),\;k=1,2,\ldots.
\]
Then the identity $w_k=1$ does not follow from the identities
\[
x^{16}=((x_1,x_2,x_3),(x_4,x_5,x_6),(x_7,x_8))=1,w_n=1,n\not=k.
\]
Here the group commutators are left normed:
\[
(u_1,\ldots,u_{n-1},u_n)=((u_1,\ldots,u_{n-1}),u_n).
\]
\end{theorem}

Another system of not finitely based group identities was found by Chander Kanta Gupta and Krasilnikov \cite{GuKr1}:

\begin{theorem}\label{groups Gupta-Krasilnikov}
The variety of groups defined by the identities
\[
(((x_1,x_2,x_3),(x_4,x_5,x_6)),x_7)=1,
\]
\[
w_m=(x^8,(y_1,y_2),\ldots,(y_{4m-3},y_{4m-2}),x^8)=1,\,m\geq 1,
\]
is not finitely based. Every identity $w_m=1$ does not follow from the other identities defining the variety.
\end{theorem}

Theorems \ref{group example of Vaughan-Lee}  and \ref{groups Gupta-Krasilnikov} show that there are commutator identities of length 7 and 8
which are not Spechtian. To the best of our knowledge the following is still an open problem, see the comments in \cite{V-L2}.

\begin{problem}
Does the identity $((x_1,x_2,x_3),(x_4,x_5,x_6))=1$ satisfy the Specht property?
\end{problem}

Kozhevnikov \cite{Kozh1} found an independent system which is simpler than the system of Adyan in Theorem \ref{example of Adyan}:

\begin{theorem}\label{system of Kozhevnikov}
The group identities $(x^p,y^p)^n=1$ form an independent system of identities for sufficiently large odd $n$
(for example, $n>10^{10}$) equal to the power of a prime number. Here $p$ runs through the set of prime numbers.
\end{theorem}

The simplest example of infinitely based varieties of groups belongs to Kleiman \cite{Klei1} and Bryant \cite{Br}.

\begin{theorem}\label{Kleiman-Bryant}
The variety ${\mathfrak B}_4{\mathfrak A}_2$ is not finitely based.
\end{theorem}

Recall that the group $G$ belongs to ${\mathfrak B}_4{\mathfrak A}_2$ if it has a normal subgroup $H$ which is abelian of exponent 2
and the factor group $G/H$ is of exponent 4. The variety has a basis of identities
\[
f_n=(x_1^2\cdots x_n^2)^4=1,\quad n=1,2,\ldots.
\]
In his paper Bryant proves a more general result: the variety ${\mathfrak B}_4{\mathfrak B}_{2k}$, $k$ odd, and many of its subvarieties are also infinitely based.

We shall complete this section with the ``continual hedgehog'' of Kleiman \cite{Klei2}, \cite{Klei3}.

\begin{theorem}\label{continual hedgehog' of Kleiman}
There exist continually many pairwise different varieties of groups ${\mathfrak V}_{\alpha}$, $\alpha\in\mathbb R$, with the same intersection
${\mathfrak V}_{\alpha}\cap {\mathfrak V}_{\beta}\subset {\mathfrak A}$ for all pairs $(\alpha,\beta)$, $\alpha\not=\beta$,
where $\mathfrak A$ is the abelian variety of groups.
\end{theorem}

\section{Infinitely based varieties of Lie algebras}

\subsection{Lie algebras in characteristic 2}

The first example of a variety of Lie algebras without a finite basis of its polynomial idenetities was given by Vaughan-Lee \cite{V-L3} in 1970.

\begin{theorem}\label{center-by-metabelian-1970}
The variety of Lie algebras over a filed of characteristic $2$
defined by the identities
\[
[[x_1,x_2],[x_3,x_4],x_5]=0,
\]
\[
[[x_1,x_2,x_3,\ldots,x_n],[x_1,x_2]]=0,\quad n=3,4,\ldots,
\]
is not finitely based. (The Lie commutators are left normed, e.g. $[u,v,w]=[[u,v],w]$.)
\end{theorem}

Volichenko \cite{Vol1} added to the system of identities in Theorem \ref{center-by-metabelian-1970} one more identity
to obtain an infinitely based locally finite variety.
(Locally finite = every finitely generated algebra is finite-dimensional.)

In the same paper \cite{V-L3} Vaughan-Lee gave an example of a finite-dimensional Lie algebra without a finite basis of its polynomial identities.
It is remarkable, that this is the Lie algebra $M_2^{(-)}(K)$ of $2\times 2$ matrices.

\begin{theorem}\label{Vaughan-Lee - M2}
If the field $K$ of characteristic $2$ is infinite, then the Lie algebra $M_2(K)^{(-)}$ of $2\times 2$ matrices
(with operation $[u,v]=uv-vu$) has no finite basis of its polynomial identities.
\end{theorem}

A basis of the identities of $M_2(K)^{(-)}$ over an infinite field of characteristic 2 was found in the Ph.D. Thesis \cite{Dr2} of the author of the present paper
and independently by Lopatin \cite{Lop}.

\begin{theorem}\label{basis of PI of M2, char 2}
The Lie algebra $M_2(K)^{(-)}$, ${\rm char}(K)=2$, $\vert K\vert=\infty$, has the following infinite basis of polynomial identities:
\[
[[x_1,x_2],[x_3,x_4],x_5]=0,
\]
\[
[[x_1,x_2,x_3,\ldots,x_n],[x_1,x_2]]=0,\quad n=3,4,\ldots,
\]
\[
[[x_1,x_2,x_5,\ldots,x_n],[x_3,x_4]]+[[x_1,x_3,x_5,\ldots,x_n],[x_4,x_2]]
\]
\[
+[[x_1,x_4,x_5,\ldots,x_n],[x_2,x_3]]=0,\quad n=4,5,\ldots.
\]
\end{theorem}

In 1975 Vaughan-Lee \cite{V-L5} gave another example of an infinitely based variety of Lie algebras.

\begin{theorem}\label{Vaughan-Lee 1975}
 If ${\rm char}(K)=2$, then the variety of Lie algebras
defined by the identities
\[
[[x_1,x_2,x_3],[x_4,x_5,x_6]]=0,
\]
\[
[x_{n+1},x_{n+2},x_{n+3},[x_1,x_2],[x_2,x_3],\ldots,[x_{n-1},x_n],[x_n,x_1]]=0,
\]
$n=2,3,\ldots$, is not finitely based.
\end{theorem}

Although not explicitly stated, the key moment of the proof that in characteristic 2 the variety of Lie algebras
${\mathfrak A}{\mathfrak N}_2$ defined by the identity
\[
[[x_1,x_2,x_3],[x_4,x_5,x_6]]=0
\]
is not Spechtian is that the system of weak polynomial identities
\[
[x,y,z]=0,[x_1,x_2][x_2,x_3]\cdots[x_{n-1},x_n][x_n,x_1]=0,n=2,3,\ldots,
\]
is not finitely based.

Over a field of characteristic different from 2 the three-dimensional Lie algebra $sl_2(K)$ of the traceless $2\times 2$ matrices is simple.
But in characteristic 2 the diagonal of the matrices in $sl_2(K)$ is a scalar matrix. Hence $sl_2(K)$ is nilpotent.
There is a simple tree-dimensional Lie algebra $(K^3,\times)$ similar to the three-dimensional real  Lie algebra with the vector product.

\begin{problem}\label{3-dim infinite basis}
If ${\rm char}(K)=2$, $\vert K\vert=\infty$,
is the basis of the polynomial identities of Lie algebra $(K^3,\times)$ finite?
\end{problem}

\subsection{Lie algebras in characteristic $p>2$}

The results of Vaughan-Lee \cite{V-L3} in characteristic 2 were transferred by the author \cite{Dr1} and Kleiman (unpublished)
to Lie algebras in characteristic $p>2$.

\begin{theorem}\label{Drensky-Kleiman}
If ${\rm char}(K)=p>0$, then the variety defined by
\[
[[[x_1,x_2],[x_3,x_4]],[[x_5,x_6],[x_7,x_8]]]=0,
\]
\[
[[x_1,x_2],\ldots,[x_{2p-1},x_{2p}],x_{2p+1}]=0,
\]
\[
[x_1,x_2,x_3,\ldots,x_n]({\rm ad}[x_1,x_2])^{p-1}=0,\quad n=3,4,\ldots,
\]
is not finitely based.
(Here $u({\rm ad}v)=[u,v]$.)
\end{theorem}

\begin{theorem}\label{locally finite char p} {\rm (Drensky \cite{Dr4})}
If we add to the system of identities from the previous theorem the identity
\[
x_1(\rm{ad}x_2)^{p^2+2}=0
\]
we shall obtain an infinitely based locally finite variety.
\end{theorem}

Over an infinite field of characteristic $p>2$ the author \cite{Dr1} obtained an analogue of Theorem \ref{Vaughan-Lee - M2}.

\begin{theorem}\label{Lie dim<infty, char p}
If ${\rm char}(K)=p>0$ and $\vert K\vert=\infty$, then there exists a finite-dimensional Lie algebra over $K$
with an infinite basis of its polynomial identities.
\end{theorem}

\begin{problem}
Are the varieties of Lie algebras over a field of characteristic $0$ finitely based?
\end{problem}

\subsection{Group graded Lie algebras}

For any abelian group $G$ an ordinary Lie algebra $L$ can be considered as a $G$-graded Lie algebra which has only a neutral component, i.e. $L=L_0$.
Hence any counterexample to the Specht problem for ordinary Lie algebras can be considered as a counterexample for $G$-graded Lie algebras.
We already mentioned that the papers by
Martinez Correa and Koshlukov \cite{MCK} and Martinez Correa and Yasumura \cite{MCY} contain positive results for graded Lie algebras with nontrivial grading.
They also obtain counterexamples in the case of characteristic 2.

\begin{theorem}\label{char 2 UT3} Let $K$ be an infinite field and let ${\rm char}(K)=2$.
Then the variety of $G$-graded Lie algebras generated by the Lie algebra $UT_3^{(-)}(K)$
does not satisfies the Specht property, where:

{\rm (i) (Martinez Correa and Koshlukov \cite{MCK})} $G={\mathbb Z}_3$ and $UT_3^{(-)}(K)$ is canonically ${\mathbb Z}_3$-graded.

{\rm (ii) (Martinez Correa and Yasumura \cite{MCY})} $G={\mathbb Z}_2$ and $UT_3^{(-)}(K)$ is
equipped with the elementary ${\mathbb Z}_2$-grading induced by $\deg(e_{12})=\deg(e_{23})=1$.
\end{theorem}

In both cases the infinitely based system of identities which defines the infinitely based subvariety consists of the identities
\[
u_n=[y^{(1)},x_1^{(0)},\ldots,x_n^{(0)},y^{(1)}]=0,\;n=1,2,\ldots.
\]

Finally, Drensky, Koshlukov and Martinez Correa \cite{DrKMC} constructed counterexamples to the Specht problem for ${\rm char}(K)=p>2$ in the spirit of the counterexamples in \cite{Dr1}.

\begin{theorem}\label{Vesselin-Plamen-Daniela}
Let $K$ be a field of positive characteristic $p>2$ and let
${\mathfrak T}_c$ be the variety of  ${\mathbb Z}_c$-graded Lie algebras defined by the $\mathbb{Z}_c$-graded identities
\[
[x_1^{(0)}, x_2^{(0)}]=[x_1^{(i)}, x_2^{(j)}]=0,\;{\rm if}\; i+j\geq c
\]
(When the field is infinite, these identities define the variety generated
by the Lie algebra $UT_c(K)^{(-)}$ with its canonical ${\mathbb Z}_c$-grading.)

{\rm (i)} The subvariety ${\mathfrak B}$ of the variety ${\mathfrak T}_{p+1}$ defined by the ${\mathbb Z}_{p+1}$-graded polynomial identities
\[
f_{n}(x_0^{(0)},\ldots, x_n^{(0)}, y^{(1)})=[y^{(1)},x_0^{(0)},\ldots, x_n^{(0)}]{\rm ad}^{p-1} y^{(1)}=0,\; n=0,1,2,\ldots,
\]
does not have a finite basis of its polynomial identities.

{\rm (ii)} The subvariety of ${\mathfrak B}$  defined by the polynomial identities
\[
y^{(i)}{\rm ad}^p x^{(0)}=0,\,\, i=1,\ldots, p,
\]
is locally finite and does not have a finite basis of its polynomial identities.

{\rm (iii)} If the field $K$ is infinite, then the variety $\mathfrak{T}_{p+1}$ contains a Lie algebra of dimension $2p+3$ over $K$
whose ${\mathbb Z}_{p+1}$-graded identities are not equivalent to a finite system of ${\mathbb Z}_{p+1}$-graded identities.
\end{theorem}

\section{Infinitely based varieties of associative algebras}

The first examples of infinitely based varieties of associative algebras in positive characteristic were presented in 1999 by Belov \cite{Bel1} (see also \cite{Bel2})
and almost in the same time by Grishin in characteristic 2 \cite{Gr2} (see also \cite{Gr1} and \cite{Gr3}) and Shchigolev \cite{Shch}.

\begin{theorem}\label{example of Belov} {\rm (Belov \cite{Bel1})}
Let $K$ be a field of characteristic $p>0$ and let $q=p^s$, $s>1$.
Define the polynomial identities
\[
R_n=[[u,t],t]\prod_{i=1}^nx_i^2
([t,[t,v]][[u,t],t])^{q-1}[t,[t,v]]=0
\]
if $p=2$ and
\[
R_n=[[u,t],t]\prod_{i=1}^n\left(x_i^{p-1}y_i^{p-1}[x_i,y_i]\right)
([t,[t,v]][[u,t],t])^{q-1}[t,[t,v]]=0
\]
for $p>2$.
Then the system of polynomial identities $R_n=0$, $n=1,2,\ldots$, is not finitely based.
\end{theorem}

\begin{theorem}\label{example of Grishin} {\rm (Grishin \cite{Gr2})}
Let $K$ be a field of characteristic $2$. Then the variety of associative algebras
defined by the polynomial identities
\[
g_0=y_1^4z_1^4z_2^4y_2^4y_1^4z_1^4z_2^4y_2^4=0,
\]
\[
g_n=y_1^4z_1^4x_1^2\cdots x_n^2z_2^4y_2^4y_1^4z_1^4x_{n+1}^2\cdots x_{2n}^2z_2^4y_2^4=0,\; n=1,2,\ldots,
\]
is not finitely based.
\end{theorem}

\begin{theorem}\label{example of Shchigolev} {\rm (Shchigolev \cite{Shch})}
Let $K$ be a field of characteristic $p>0$, let $k,l\geq 2p^3$ be some integers and let
\[
Q_n=z_1^2\dots z_n^2,\; {\rm if}\;p=2,
\]
\[
Q_n=(z_1^{p-1}z_2z_1z_2^{p-1})\cdots(z_{2n-1}^{p-1}z_{2n}z_{2n-1}z_{2n}^{p-1},\;{\rm if}\;p>2,
\]
\[
f_n=\prod_{j=0}^{p-1}(y_1^{p^2}x_1\cdots x_kQ_{nj}x_{k+1}\cdots x_{k+l}y_2^{p^2}),
\]
where $Q_{nj}$ is obtained from $Q_n$ by shifting the variables $z_1,\ldots,z_n$ in $Q_n$ by $n$ for $p=2$
and the variables $z_1,\ldots,z_{2n}$ in $Q_n$ by $2n$ for $p>2$.
Then the variety defined by the system of polynomial identities $f_n=0$, $n=1,2,\ldots$, does not satisfy the finite basis property.
\end{theorem}

It is interesting to mention that the examples of Grishin and Shchigolev consist of monomials and the example of Grishin consists of nil algebras of index 32, i.e.
the variety satisfies the identity $x^{32}=0$.

Maybe the simplest examples of not finitely based varieties of associative algebras are due to Gupta and Krasilnikov \cite{GuKr2} and \cite{GuKr3} in characteristic 2.

\begin{theorem}\label{two papers Gupta-Krasilnikov}
The following two systems of polynomial identities over a field of characteristic $2$ are not Spechtian:
\[
x^6=0, w_n=u_1x_1^2\cdots x_n^2u_2u_1u_2=0,  u_i=[[y_i,z_i],t_i], n\geq 0.
\]
\[
[x,y^2]x_1^2\cdots x_n^2[x,y^2]^3=0,n\geq 0.
\]
\end{theorem}

We finish this section with the following conjecture.

\begin{conjecture}
Over an infinite field $K$ of characteristic $2$ the associative algebra $M_2(K)$
of $2\times 2$ matrices does not have a finite basis of its polynomial identities.
\end{conjecture}

\section{Algebras that are nearly associative}

This is the common name of nonassociative algebras which share many properties of the associative algebras.
Typical examples are the Jordan and the alternative algebras.

An algebra $R$ is right-alternative if it satisfies the identity
$(xy)y=x(yy)$, i.e. $(x,y,y)=0$. Belkin \cite{Belk} constructed a right-alternative algebra
over an arbitrary field which does not have a finite basis of its polynomial identities.

The first example of a non-Spechtian variety of alternative algebras over a field of characteristic 2 was given by Medvedev \cite{Med2}.
His variety is defined by the identities
\[
((x_1x_2)(x_3x_4))x_5=x_1((x_2x_3)(x_4x_5))=0,
\]
\[
(\cdots((x^2y_1)y_2)\cdots y_n)x=0,\,n=4k+2\geq 6.
\]
It is interesting to mention that over a field of characteristic different from 2 this variety is finitely based.

Pchelintsev \cite{Pch2} showed that over an infinite field of characteristic 3 the variety of alternative algebras defined by the identities
\[
[[x,y],z]=((x_1x_2)(x_3x_4))x_5=x_1((x_2x_3)(x_4x_5))=0,
\]
\[
(\cdots([x_1,x_2]x_3)\cdots x_{n-2})[x_{n-1},x_n]=0,
\]
is not Spechtian. As usually, here $[x,y]=xy-yx$ is the Lie commutator.
This variety is limit, i.e. its proper subvarieties are finitely based.

\section{Finite rings and finite-dimensional algebras}

In 1976 Polin \cite{Poli} constructed a finite linear algebra without a finite basis of its identities.

\begin{theorem}\label{the algebra of Polin}
Over every finite field $K$ there is a finite-dimensional algebra $R$ which is left nilpotent, i.e. satisfies
$x_1(x_2x_3)=0$, such that ${\rm var}(R)$ is not finitely based.
\end{theorem}

In 1978 Lvov \cite{Lv2} suggested the following construction which in the sequel was used and improved
to construct finite-dimensional nonassociative algebras without finite bases of their polynomial identities.

For any field $K$ the algebra of Lvov is of the form
\[
R=V_2\oplus M_2(K),
\]
where $M_2(K)$ acts from the right on the two-dimensional vector space $V_2$ and the multiplication in $R$ is defined by
\[
(v_1+a_1)(v_2+a_2)=v_1a_2,\quad m_1,m_2\in V,a_1,a_2\in M_2(K).
\]
Again, $R$ satisfies the polynomial identity $x(yz)=0$.

More generally, Lvov suggested to consider the algebra $R=V_d\oplus M_d(K)$, where $V_d$ is a $d$-dimensional vector space
with a similar multiplication rule.

\begin{theorem}\label{the method of Lvov}
Over an arbitrary field $K$ the algebra $R=V_d\oplus M_d(K)$ has a finite basis of its polynomial identities if and only if
the T-ideal $T(M_d(K))\triangleleft K\langle X\rangle$ is generated as an ordinary ideal (not as a T-ideal) by polynomials of bounded degree.
\end{theorem}

Since this does not hold for the $2\times 2$ matrix algebra $M_2(K)$ over any field one immediately obtains:

\begin{corollary}\label{the algebra of Lvov}
The polynomial identities of the six-dimensional nonassociative algebra $R=V_2\oplus M_2(K)$ over any field $K$ do not follow from a finite number of identities.
\end{corollary}

Over a field $K$ of characteristic 0
Maltsev and Parfenov \cite{MaltsPar} constructed a 5-dimensional algebra with an infinite basis of its polynomial identities.

\begin{theorem}\label{Maltsev and Parfenov}
Let $R$ be the $5$-dimensional algebra with basis $\{e_1,e_2,e_3,e_4,e_5\}$ over a field $K$ of characteristic $0$
with nonzero multiplications between the  basis elements of the algebra
\[
e_1e_3=e_2e_3=e_1,e_1e_4=e_2e_4=e_2,e_1e_5=e_1,e_2e_5=e_2.
\]
Then the identities
\[
x(yz)=0,x[y,z]x_1\cdots x_n[u,v]=0,n=0,1,2,\ldots,
\]
form a minimal basis of the polynomial identities of $R$. Here
the products are left normed and $a[b,c]$ is a shortcut for $abc-acb$.
\end{theorem}

It is easy to see that translated in the construction of Lvov (and changing the basis of $R$ as a vector space) this is the algebra $R=V_2\oplus UT_2(K)$, where
$UT_2(K)$ is the algebra of $2\times 2$ upper triangular matrices.

The basis of the identities of the algebra of Maltsev and Parfenov over a finite field was given by Isaev \cite{Isa2}
and over an arbitrary infinite field of positive characteristic by Isaev and Kislitsin \cite{IsKi1}.

Finally, a four-dimensional algebra without a finite basis of its identities was constructed by Isaev and Kislitsin
(in \cite{IsKi2} over finite fields and in \cite{IsKi3} over infinite fields).

\begin{theorem}\label{Isaev-Kislitsin, dim=4}
The following $4$-dimensional algebra has an infinite basis of its identities over any field
\[
R=V_2+{\rm span}\{e_{11}+e_{12},e_{22}\}, \dim(V_2)=2.
\]
\end{theorem}

A similar construction was used by Isaev \cite{Isa3} who found two Spechtian varieties ${\mathfrak M}_2$ and ${\mathfrak M}_3$
over a field of characteristic 0
such that their union ${\mathfrak M}_2\cup{\mathfrak M}_3$ is not Spechtian.
Let $R=V_2\oplus M_2(K)$ be the algebra of Lvov in Corollary \ref{the algebra of Lvov}. Isaev consider the subalgebras $R_1$ and $R_2$ of $R$ and the subalgebra $R_3$ of $R_1$:
\[
R_1=V_2+{\rm span}\{e_{11},e_{21}\},\,R_2=V_2+{\rm span}\{e_{11},e_{12}\},\,R_3={\rm span}\{v_1,v_2+e_{11},e_{21}\}.
\]
Kislitsin \cite{Ki1} found the bases of the identities of the representations of the subalgebras $\{e_{11},e_{21}\}$ and $\{e_{11},e_{12}\}$ of ${\rm End}(V_2)=M_2(K)$.
As a consequence, the bases of the identities of $R_1$ and $R_2$ consist, respectively, of
\[
x_1(x_2x_3)=x_1x_2[x_3,x_4]=0\;{\rm and}\;
x_1(x_2x_3)=x_1[x_2,x_3]x_4=0.
\]
A basis of the identities of $R_3$ consists of
\[
x_1(x_2x_3)=x_1[x_2,x_3]x_4 - x_4[x_2,x_3]x_1 = 0,
\]
\[
s_3 = x_1[x_2,x_3] + x_2[x_3,x_1] + x_3[x_1,x_2] = 0.
\]
(As in Theorem \ref{Maltsev and Parfenov} the products are left normed.)
Both varieties ${\mathfrak M}_2={\rm var}(R_2)$ and ${\mathfrak M}_3={\rm var}(R_3)$ satisfy the Specht property but their union ${\mathfrak M}_2\cup{\mathfrak M}_3$
generated by $R_2\oplus R_3$ has a basis of identities
\[
x_1(x_2x_3)=x_1[x_2,x_3]x_4 - x_4[x_2,x_3]x_1 = 0,
\]
\[
x_1[x_2,x_3]x_4 + x_1[x_3,x_4]x_2 + x_1[x_4,x_2]x_3 = 0,\;
x_1x_2[x_3,x_4]x_5 = 0,
\]
\[
x_1[x_2,x_3]y_1y_2\cdots y_n[x_4,x_5] = 0,\;n = 0, 1, 2,\ldots .
\]
which do not follow from a finite subsystem.

Similar results are obtained also for classes of nearly associative algebras.
Over an arbitrary field Isaev \cite{Isa1} gave an example of a finite-dimensional right-alternative algebra without a finite basis of its polynomial identities.
Again, his construction is inspired by the construction of Lvov \cite{Lv2}. Using the constructions and methods in \cite{Lv2} and \cite{Isa1}, over any field
Ismailov and Umirbaev \cite{IsUm} constructed a finite-dimensional right-symmetric algebra without a finite basis of its identities.

The existence of infinitely based varieties generated by algebras of dimension 6, 5 and 4 in Corollary \ref{the algebra of Lvov} and Theorems \ref{Maltsev and Parfenov}
and \ref{Isaev-Kislitsin, dim=4} naturally inspirits the following problem which generalizes Problem \ref{3-dim infinite basis}.

\begin{problem}\label{3-dim nonassociative infinite basis}
Does there exist a three-dimensional nonassociative algebra $R$ without a finite basis of its polynomial identities?
\end{problem}

One possible candidate for such an algebra is the three-dimensional algebra with basis $\{e_0,e_1,e_2\}$ and multiplication between the basis elements
defined as the multiplication of the three-element binary magma $M$ of Murskii in Theorem \ref{Murski n=3}.

The proofs for the infinitely based varieties generated by algebras of dimension 6,5 and 4 are based essentially on the infinite basis property of weak polynomial identities.
We shall finish this section with the following problem.

\begin{problem}\label{weak 3-commutator}
Is the weak polynomial identity $[x,y,z]=0$ Spechtian in characteristic $p>2$?
\end{problem}

The negative answer would give an example of infinitely based subvariety
of the variety of Lie algebras ${\mathfrak A}{\mathfrak N}_2$ in characteristic $p>2$
which would be a nice complement to the result by Vaughan-Lee in Theorem \ref{Vaughan-Lee 1975} in characteristic 2 and by Volichenko in Theorem \ref{Volichenko week PI} in characteristic 0.

\section{Limit varieties}

The variety of algebraic systems is limit (or just-non-finitely based) if it is not finitely based
and all of its proper subvarieties are finitely based.
By the Zorn lemma every infinitely based variety has a limit subvariety.

We shall survey some results on explicit examples of limit varieties.
The examples do not depend on the Zorn lemma.

\subsection{Groups}

The first examples of infinitely based varieties of groups guarantee that there are limit varieties of groups.
It follows from Theorem \ref{group example of Vaughan-Lee} of Vaughan-Lee that there exists a limit variety of 2-groups.
Theorem \ref{theorem of Olshanski} of Ol'shanskij and Theorem \ref{example of Adyan} by Adyan do not show that there exist other limit varieties.
Newman \cite{New} showed that for any prime $p>2$ there is an infinitely based variety of $p$-groups.

\begin{theorem}\label{example of Newman}
Let $p$ be a prime, ${\mathfrak A}_p$ be the variety of abelian groups of exponent $p$,
${\mathfrak N}_p$ be the variety of nilpotent groups of class $p$, defined by the identity $(x_1,\ldots,x_{p+1})=1$,
and ${\mathfrak T}_p$ be the variety generated by the nonabelian group of order $p^3$ and exponent $p$. Let
\[
{\mathfrak W}=({\mathfrak A}_p{\mathfrak A}_p\cap{\mathfrak N}_p){\mathfrak T}_p.
\]
The subvariety $\mathfrak V$ of $\mathfrak W$ defined by the identities
\[
v_k=((x,y,z),(x,y,z)^{u_k},\ldots,(x,y,z)^{u_k^{p-1}})=1,k=1,2,\ldots,
\]
where $u_k=(x_1,x_2)\cdots(x_{2k-1},x_{2k})$ is not finitely based. The identities $v_k=1$ are an independent system
if we consider $\mathfrak V$ as a subvariety of $\mathfrak W$.
\end{theorem}

\begin{corollary}
For any prime $p$ there is a limit variety of $p$-groups. Hence there are infinitely many limit varieties.
\end{corollary}

Theorem \ref{example of Newman} was improved by Kozhevnikov \cite{Kozh2}, \cite{Kozh3},
who established that there are uncountably many limit varieties of groups.

No explicit examples are known for limit varieties of groups.

\subsection{Semigroups}

The first example of a limit variety of semigroups was constructed by Volkov \cite{Volk1}.

\begin{theorem}\label{limit variety of Volkov}
Let $\mathfrak W$ be the variety of semigroups defined by the identities
\[
xyzw=xzyw,x^2y^2=y^2x^2,x^2zy^2=y^2zx^2
\]
and let $\mathfrak V$ be the subvariety of $\mathfrak W$ defined by the identities
\[
yx^2y=y^2x^p,\quad p\geq 2\;{\rm prime}.
\]
Then $\mathfrak V$ is the only limit subvariety of $\mathfrak W$.
\end{theorem}

More examples of infinite series of limit varieties (some of which are generated by finite semigroups) can be found in the papers by
Poll\'ak \cite{Poll}, M.V. Sapir \cite{MVSap}, Lee and Volkov \cite{LeeVol}.
Klejman \cite{Klej} showed that the 6-element Brandt semigroup $B_2^1$ generates a
limit variety of inverse semigroups and this is the only limit variety of inverse semigroups which is not a limit variety of groups.

The hunting for limit varieties of monoids continues also nowadays.
See for example the following papers and the references there:
Lee \cite{Lee1}, Zhang and Luo \cite{ZhL}, Gusev \cite{Gus1},
Olga B. Sapir \cite{OSap1} and \cite{OSap2}, Gusev  \cite{Gus2},
Gusev and O.B. Sapir \cite{GusSap}.
The recent book \cite{Lee4} by Lee contains a lot of information on the finite basis problem for semigroups and related topics.
The companion website \cite{AACLR2} of the survey \cite{AACLR1} by
J. Ara\'ujo, J.P. Ara\'ujo, Cameron, Lee, and Raminhos contains the bases of the identities of all varieties generated by semigroups of order $\leq 4$.

A variety of algebraic systems is irredundantly based if it can be defined by an independent system of identities.
For example, the defining system of identities
\[
f_n=f(x_1,\ldots,x_n)=(x_1^2\cdots x_n^2)^4=1,\quad n=1,2,\ldots,
\]
of the variety of groups ${\mathfrak B}_4{\mathfrak A}_2$ is not independent because each identity $f_n=1$
implies all $f_m=1$ for $m<n$.

Similar examples are known also for semigroups, see Mashevitskij \cite{Mash},
M.V. Sapir \cite{MVSap}, Lee \cite{Lee2}.

The paper by Lee \cite{Lee2} shows that the variety generated by the semigroup
\[
L=\langle e,f\mid e^2=e,f^2=f,efe=0\rangle=\{0,e,f,ef,fe,fef\}
\]
does not have an independent system of defining identities.
Additionally, the same property holds for the variety of semigroups generated by the direct product $L\times {\mathbb Z}_n$, where
${\mathbb Z}_n$ is the cyclic group of order $n$.
The variety generated by $L$ is minimal in the following sense. All semigroups of order $\leq 5$ are finitely based, there are only four non-finitely based semigroups of order 6
and $L$ is one of them.
The semigroup $L$ can also be linearly represented by $3\times 3$ upper triangular binary matrices.
Other interesting properties of this semigroup can be found in \cite{Lee9}.

\subsection{Lie algebras}

Due to Volichenko \cite{Vol1} only one example of a limit variety of Lie algebras in characteristic 2 is explicitly known.
It is a subvariety of the center-by-abelian variety $[{\mathfrak A}^2,{\mathfrak E}]$ defined by the identity
\[
[[x_1,x_2],[x_3,x_4],x_5]=0.
\]

\subsection{Finite magmas}
The description of the lattice of subvarieties of the variety generated by the magma of Lyndon given in Theorem \ref{paper Lee5}
implies that this variety is limit.

\subsection{Finite-dimensional algebras in characteristic 0}

In the paper by the author \cite{Dr3} a limit variety of anticommutative algebras was constructed.

\begin{theorem}\label{anticommutative limit}
Over a field of characteristic $0$ the variety of algebras $\mathfrak V$ defined by the following identities is limit:
\[
xy+yx=0,(xy)(zt)=0,x_1x_2y_1y_2x_3-x_1x_2y_2y_1x_3=0,
\]
\[
s_3(x_1,x_2,x_3)x_1^{n-5}[x_1,x_2]=0,n=5,6,\ldots,
\]
\[
\sum_{\sigma\in S_4}{\rm sign}(\sigma)x_{\sigma(1)}x_{\sigma(2)}x_{\sigma(3)}x_1^{n-4}x_{\sigma(4)}=0,n=4,5,\ldots.
\]
\end{theorem}

Here the products are left normed: $xyz=(xy)z$,
$x[y,z]=xyz-xzy$ and $s_3(x_1,x_2,x_3)$ is the standard polynomial of degree 3.
The proof uses representation theory of the general linear group.
The variety $\mathfrak V$ from the previous theorem is generated by the anticommutative six-dimensional algebra with a vector space basis $\{a_1,a_2,a_3,b,c,g\}$
and nonzero products of the basis elements
\[
a_1a_2=b,ba_3=c,bg=b.
\]

The proof of the following theorem is similar to the proof in the anticommutative case in Theorem \ref{anticommutative limit}.

\begin{theorem}
The variety $\mathfrak V$ defined by the identities
\[
xy-yx=0,(xy)(zt)=0,x_1x_2y_1y_2x_3-x_1x_2y_2y_1x_3=0,
\]
\[
\sum_{\sigma\in S_2}{\rm sign}(\sigma)x_{\sigma(1)}[x_1,x_2]x_{\sigma(2)}x_1^{n-6}[x_1,x_2]=0,n=6,7,\ldots,
\]
\[
\sum_{\sigma\in S_3}{\rm sign}(\sigma)x_{\sigma(1)}[x_1,x_2]x_1^{n-5}x_{\sigma(2)}x_{\sigma(3)}=0,n=5,6,\ldots.
\]
is limit and is generated by the commutative six-dimensional algebra with a vector space basis $\{a_1,a_2,a_3,b,c,g\}$
and nonzero products of the basis elements
\[
a_1a_2=b,ba_3=c,bg=b.
\]
\end{theorem}

\subsection{Weak polynomial identities}
Very recently Kislitsin \cite{Ki2} constructed an example of a limit variety of pairs which is almost commutative and almost Engel,
i.e. all of its proper subvarieties satisfy the weak commutativity identity and the weak Engel identity,
but the variety itself does not satisfy them.

\section{The Specht problem for other algebraic constructions}

Up till now we considered the Specht problem for different kinds of algebraic systems.
There are also other algebraic constructions for which one can study the Specht problem.
Such mathematical objects are:

$\bullet$ {\bf Varieties of group representations.} This topic is intensively studied by the school of Boris Plotkin in Riga,
see the survey by Plotkin \cite{Plo}, the books by Plotkin and Vovsi \cite{PloVo} and Vovsi \cite{Vov}.
The topic is also related with problems in the theory of varieties of groups, see e.g. Gupta and Krasilnikov \cite{GuKr1}.

$\bullet$ {\bf Varieties of representations of Lie algebras.} Again, the school of Plotkin was very active in this direction,
see the books \cite{PloVo} and \cite{Vov}. The paper by Krasilnikov and Shmelkin \cite{KrShm}
studies varieties of representations of Lie algebras over a field of characteristic 0 and over an infinite field of characteristic 2.
The obtained results are in the spirit of the results in Krasilnikov \cite{Kra} and Vaughan-Lee \cite{V-L3} but the proofs are based on different ideas.

$\bullet$ {\bf Endomorphism semirings of semilattices.}
For details we refer to a series of papers by Crvenkovi\'c and Dolinka \cite{CrDol}, Dolinka \cite{Dol1}, \cite{Dol2} and \cite{Dol3} and Dolinka, Gusev and Volkov \cite{DGV}.

\section{Further readings (on counterexamples\\ to the Specht problem only)}

Here are some books and surveys (listed chronologically) which contain results on the Specht problem for different algebraic systems
(mainly groups and Lie and associative rings and algebras):
Hanna Neumann \cite{HNeu}, Bahturin \cite{Ba1}, \cite{Ba2}, \cite{Ba3}, Bahturin and Olshanskii \cite{BaOl2},
Drensky \cite{Dr5}, Kanel-Belov and Rowen \cite{K-BR},
Belov-Kanel, Rowen, and Vishne \cite{B-KRV},
Kanel-Belov, Karasik, and Rowen \cite{K-BKR},
Aljadeff, Giambruno, Procesi, and Regev \cite{AGPR}.

\end{document}